\documentclass{article}
\usepackage{amsfonts}
\usepackage{amsmath}
\usepackage{amssymb}
\usepackage{color}
\usepackage{graphicx}
\newcounter{excnt}
\newtheorem{thm}{Theorem}

     {\noindent {\bf Proof:}}{\nopagebreak\begin{flushright}{\bf q.e.d.}\end{flushright}}
\newenvironment{example}%
     {\stepcounter{excnt}
      \noindent {\bf Example \arabic{excnt}} \phantom{i}}{\vspace{0.4cm}}

\title{Desingularization Algorithms\\
       A Comparison from
       the Practical Point of View}
\author{Rocio Blanco (Universidad de Castilla-La Mancha. \\ Depto. de Matem\'aticas. Facultad de Educaci\'on. Cuenca, Spain)\\ \\
        Anne Fr\"uhbis--Kr\"uger (Institut f. Alg. Geometrie, \\
                                  Leibniz Universt\"at Hannover, Germany)}

\begin{document}

\maketitle

\begin{abstract}
          Over the last decade, implementations of several desingularization
          algorithms have appeared in various contexts. These differ as 
          widely in their methods and in their practical efficiency as 
          they differ in the situations in which they may be applied. 
          The purpose of this article is to give a brief overview
          over a selection of these approaches and the applicability of 
          the respective implementations in {\sc Singular}.
\end{abstract}

\section{Introduction}

The problem of resolution of singularities has stimulated mathematical
research for more than a century now, ranging from first results for
curves in the late 19th century through to the 
case of surfaces (e.g. \cite{Jung}, \cite{Lipman}) and 3-folds 
(e.g. \cite{Zariski}), culminating in the famous proof of Hironaka 
for the general case in characteristic zero in 1964 \cite{Hir}. But this 
result did not end the interest in the problem: First of all, the case of 
positive characteristic is still an open problem and an active field of 
research; on the other hand, the fact that Hironaka's proof was highly 
non-constructive, led to the question of the existence of algorithmic 
approaches in characteristic zero and to the development thereof 
(e.g. \cite{BM1},\cite{Vil},\cite{EH} etc.). \\

In its simplest form the problem of desingularization can be phrased as:
Given a variety $X$, determine a non-singular variety $\tilde{X}$ and 
a proper birational morphism $\pi: \tilde{X} \longrightarrow X$ which 
leaves $X \setminus Sing(X)$ unchanged. In general, it is a rather
complicated task to construct such a morphism $\pi$ by an iteration of
improvements of $X$ through blow-ups in appropriate centers,
where the choice of center turns out to be the crucial point. In certain
special situations, however, like binomial varieties or curves and surfaces,
significant simplifications in the choice of centers and the computation
of the blow-ups or  the additional use of normalization steps can lead 
to faster algorithms which sometimes also provide an $\tilde{X}$ which 
is covered by fewer charts than the one constructed by the general algorithm. 
This has led to implementations of the various types of algorithms of 
which at least one for each type is available as a {\sc Singular} library
as of version 3.2. Being offered a variety of algorithms for very similar
tasks is certainly a nice feature for the experienced user, but novice users
are easily overwhelmed by it without hints to guide their choice. Therefore
a thorough discussion of the structural differences and the resulting 
pros and cons is the topic of this article. As it is not possible to 
explain all these algorithms in detail in one article (this would rather
require a book) or at least do each of them justice in a short description,
we only highlight the main ideas and general structure of each type of 
approach and give the corresponding {\sc Singular} commands and an 
interpretation of the resulting output in explicit examples.\\
 
For comparing algorithms, computer science usually suggests
time and space complexity considerations. In this case, however, it 
has been shown in \cite{BGMW} 
% Bierstone, Grigoriev, Milman, Wlodarczyk: effective Hironaka resolution and its complexity to appear in Asian J Math. -- nachsehen!
that the complexity of the general resolution algorithm is in the 
complexity class ${\mathcal E}^{m+3}$ for a $m$-dimensional 
variety\footnote{To give the reader an orientation about Grzegorczyk 
complexity classes, the class ${\mathcal E}^1$ contains all linear functions,
${\mathcal E}^2$ all polynomial functions and ${\mathcal E}^3$ all towers
of exponential functions .}, which for a computer scientist basically 
translates to 'it is completely hopeless to calculate anything with this 
algorithm'.
On the other hand, algebraic geometers know from experience with Groebner
basis computations or primary decomposition (both in ${\mathcal E}^3$) 
that many examples, which they want to compute, are quite far from the 
maximum complexity of the algorithm. This suggests to select a set of 
typical examples from practical applications and use these for the 
comparison, which is our point of view in this article.\\

After very briefly sketching the common main ideas of all Hironaka style 
approaches and giving a {\sc Singular} example thereof in the section 
\ref{sec2}, we consider the structurally much simpler binomial and toric
case in the section \ref{sec3}. In the subsequent section \ref{sec4}, we then
consider specialized algorithms for low-dimensional varieties and then
conclude the article with a comparison of the algorithms in section \ref{sec5}. 

The authors are grateful to Santiago Encinas, Herwig Hauser, Ana Bravo, 
Orlando Villamayor, Santiago Zarzuela, Nobuki Takayama, Mohammed Barakat,
Hans Sch\"onemann and Gerhard Pfister for fruitful discussions, some 
concerning the mathematical side of this article, some related to 
software development and applicability of the respective {\sc Singular} 
libraries. 

\section{Algorithms refining Hironaka's approach in the general case} 
\label{sec2}

Hironaka's groundbreaking 1964 article \cite{Hir} on resolution of 
singularities in arbitrary dimension in characteristic zero contains both 
constructive and non-constructive arguments. Seeking better understanding 
of this proof, these non-constructive steps were then studied in detail and
filled with explicit algorithmic constructions leading to the first
algorithmic proofs of desingularization in arbitrary dimension in the late
1980s/early 1990s (see e.g. \cite{BM1} \cite{Vil}), which were subsequently 
refined further (e.g. \cite{EH}). Although these approaches may differ
significantly in various constructive steps they all follow the general
structure of the original proof.\\

To explain all variants of algorithmic proofs in the flavour of Hironaka
in detail would require at least a book, so we concentrate on giving
a brief impression of the common main ideas of all these approaches and the
resulting computational tasks -- in the hope of providing the context for 
the corresponding examples in {\sc Singular}.\\

\begin{thm}[Hironaka, 1964]
For every algebraic variety over a field of characteristic zero, a 
desingularization can be achieved by a finite sequence of blow ups
in suitable non-singular centers.   
\end{thm}

From this formulation, we already see that two different algorithmic
tasks arise here: the blow up, which is a Gr\"obner basis computation
and hence does not pose any significant problems, and the choice of 
suitable centers, which is the very heart of the construction. For
our context, the gist of the choice of center can best be formulated
by means of a resolution function whose maximal locus yields the upcoming 
center.\\

To define the resolution function, we work with the ideal defining the 
algebraic variety, and we proceed by induction on the dimension of the 
ambient space.
This resolution function, which needs to be both Zariski- and infinitesimally upper semi-continuous, has the following general structure:
$$(inv_d; inv_{d-1}; \dots )$$
where each component $inv_i$ stands for the contribution by an auxiliary
ideal in ambient dimension $i$. Each of these components is then composed
of an order of the ideal (or in top dimension in some approaches the 
Hilbert-Samuel function) and some counting function which ensures normal
crossing with the exceptional divisors. The construction of the auxiliary
ideals is quite delicate and a key step in Hironaka's inductive argument.
To illustrate the flavour of this construction without going into too 
many technicalities, let us consider a power series w.r.t. a main 
variable $z$
$$f=z^k+a_1(\underline{x}) z^{k-1}+ \dots + a_k(\underline{x}) \subset
    {\mathbb C}\{z,\underline{x}\}.$$
It has order at most $k$ and this order is attained if and only if the
order of each coefficient $a_i(\underline{x})$ is at least $i$. This can
also be phrased as 
$$ord_{\underline{0}}(f) = k \;\; \Longleftrightarrow \;\;
        ord_{\underline{0}}(\langle a_1^{k!},a_2^{\frac{k!}{2}}, \dots
                                    a_k^{(k-1)!}) \rangle \geq k!$$
where the ideal generated by these powers of the coefficients now plays the
role of the first auxiliary ideal. Of course, the choice of the main
variable, i.e. the hypersurface of maximal contact,  is crucial in this 
construction and additional difficulties here arise from the fact that 
such a hypersurface exists only locally.\footnote{For an accessible and 
concise treatment of an approach of this flavor and the precise 
definitions of the technical terms used here, we refer the reader to \cite{EH}.} \\

There are two implementations of Villamayor-style algorithms for the 
general case, both in {\sc Singular}: \cite{ImplBSch} and \cite{ImplFKP} where
the latter one is part of the standard {\sc Singular} distribution and is
well suited for use with the utilities in the additional library 
\cite{ImplFKP2} for further computations with resolution data. In this 
article, we will only give a syntax example for this second implementation:\\

\begin{example}

\begin{verbatim}

> LIB "resolve.lib";           // load the library resolve 
> ring r=0,x(1..3),dp;         // ring of char 0 with 
                               // variables x(1),x(2) and x(3) 
> ideal J=x(1)^2-x(2)*x(3)^2;  // Whitney umbrella 
> list re=resolve(J);          // compute the desingularization
> size(re);                    // the output is a list 
2                              // of 2 components 
> re[1];                       // list of final charts
                               // every chart is given by
                               // its defining ring
[1]:                           // there is only one final chart
// characteristic : 0
// number of vars : 3
//      block   1 : ordering dp
//                : names    x(2) x(3) y(1)
//      block   2 : ordering C
> re[2];                       // list of all charts
[1]:
// characteristic : 0
// number of vars : 3
//      block   1 : ordering dp
//                : names    x(1) x(2) x(3)
//      block   2 : ordering C
[2]:
// characteristic : 0
// number of vars : 3
//      block   1 : ordering dp
//                : names    x(2) x(3) y(1)
//      block   2 : ordering C
\end{verbatim}

The output is a list consisting of two lists of rings. The first component gives the final charts, and the second component is the list of all charts. To see the information which is inside a concrete ring or chart it is enough to change to the respective ring and proceed as follows: 

\begin{verbatim}
> def RR=re[2][2];                         // define a new ring
> setring RR;                              // set this ring
> showBO(BO);                              // show the information
==== Ambient Space
_[1]=0                                     // the whole space
==== Ideal of Variety:                     // transform ideal
_[1]=y(1)^2-x(2)
==== Exceptional Divisors:                 // list of exceptional divisors
[1]:
_[1]=x(3)
==== Images of variables of original ring: // morphism defining
_[1]=x(3)*y(1)                             //  the blowing up
_[2]=x(2)
_[3]=x(3)
\end{verbatim}

In addition, using the additional library \cite{ImplFKP2} it is possible to obtain a picture illustrating the hierarchy of charts. We include here a Singular syntax example, but omit the picture, as it only contains two charts in this
example. (Example 2 in the next section leads to a more meaningful picture
included as figure 1.)

\begin{verbatim}
> LIB "reszeta.lib";         // library needed for visualization
> LIB "resgraph.lib";        // library needed for visualization
> list iden0=collectDiv(re); // collect all the exceptional divisors
> ResTree(re,iden0[1]);      // creates a jpg file with the resolution tree 
\end{verbatim}
\end{example}

As this general approach is obviously quite time and space consuming, it is
worthwhile to use specialized algorithms wherever those are available.

\section{Combinatorial algorithms for the binomial case} \label{sec3}

In the setting of binomial ideals or toric ideals, the situation is 
much simpler than in the general case. This allows a restriction of the
methods to combinatorial constructions and arguments. Apart from the 
classical algorithm of star subdivision as toric resolution of 
singularities (see e.g. standard textbooks about toric geometry like 
\cite{Ful}), there are combinatorial variants of both main algorithmic 
approaches to Hironaka resolution, i.e. one in the flavour of Bierstone 
and Milman \cite{BMbinom} and one in the philosophy of Villamayor \cite{preprintBlancoPartI}, \cite{preprintBlancoII}.\\

Here we will only consider the latter one, which is available as
\cite{ImplBl} in the standard distribution of {\sc Singular}.
It is mainly a combinatorial game on the exponents of the generators 
of the ideal. Roughly speaking, the algorithm looks for a smallest 
order term in the binomials and chooses one of the variables involved in 
that term to make an induction on the dimension of the ambient space. 
Then, it looks again for the next smallest order.\\ 

This combinatorial game provides an ideal generated by hyperbolic 
equations that is locally monomial. This ideal needs to be rewritten 
as a monomial one to achieve a log-resolution, see \cite{preprintBlancoPartI} 
and \cite{preprintBlancoII} for details.\footnote{This part of the algorithm 
is not yet implemented (June 2011).} Another benefit of the purely 
combinatorial approach is the applicability of the algorithm in positive
characteristic (in contrast to the general algorithms). In a bit more detail,
the philosophy of the computations can be illustrated in the scheme below: 

\begin{small}
$$ \begin{array}{ccccc} 
 \text{\underline{ Arbitrary char }} & & \underline{\text{ Char } 0 } &  & \underline{\text{Computations}} \vspace*{0.2cm} \\
 \text{Input: Binomial ideal} & \rightarrow &  \text{coefficients, exponents} &  \put(0,5){\vector(3,-2){10}} &  \\
 & & \text{type integer} & &  \text{exponents=rational numbers}  \\

 & & & &  \put(-35,0){\vector(-3,-2){10}} \text{type number}   \\
 \text{Output: List of rings } & \leftarrow &  \text{coefficients, exponents} & & \\ \text{(affine charts)} & & \text{type integer} & & \\
\end{array}  $$   
\end{small}

\begin{example}
\begin{verbatim}
> LIB "resbinomial.lib";         // load the library resbinomial
> ring r = 2,(x(1..3)),dp;       // ring of char 2 with
                                 // variables x(1),x(2) and x(3)  
> ideal J=x(1)^2-x(2)^2*x(3)^2;
> list B=binResolution(J);       // compute the binomial resolution
> size(B);                       // the output is a list of 4 components
4  
> B[1];                          // list of final charts in char 0 
> B[2];                          // list of all charts in char 0 
> B[3];                          // list of final charts in char 2
> B[4];                          // list of all charts in char 2
\end{verbatim} 
\medskip

The first two components of the output are needed for internal computations and also for  visualization of charts, since the input for the additional library \cite{ImplFKP2} has to be in characteristic zero. Then, for visualization of the resolution tree, select the first two components.

\begin{verbatim}
> list L2=B[1],B[2];
> list iden0=collectDiv(L2);
> ResTree(L2,iden0[1]);
\end{verbatim}

\begin{figure}[h]
\begin{center}
\includegraphics[width=8cm]{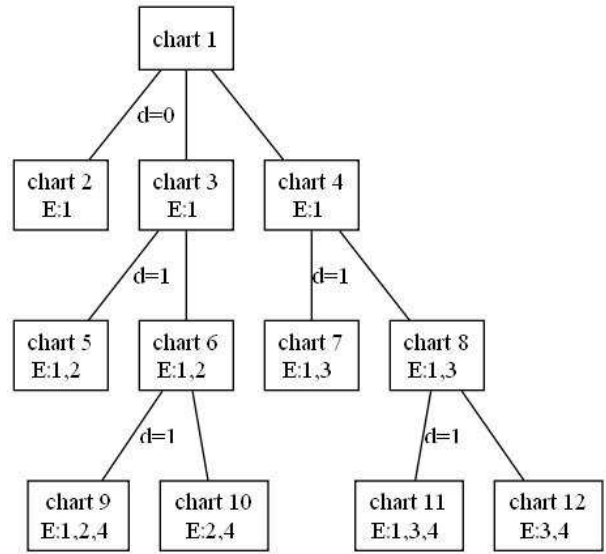}
\caption{Resolution tree of Example 2. d denotes the dimension of the center
         of the respective blow up. At each chart, $E$ lists the set of 
         exceptional divisors visible in the chart. }
\end{center}
\end{figure}
\end{example}

\section{Algorithmic resolution in low dimensions} \label{sec4}

The other main setting, in which specialized algorithms exist, is the 
situation of low dimensional varieties, i.e. curves and surfaces. As before,
these situations are significantly simpler than the general case and the
first approaches to these problems historically precede Hironaka's proof
by decades.\\
As the singular loci of curves can only be finite sets of points, 
no special precautions to choose appropriate centers need to be taken 
and a resolution can be achieved by simply using the points of the singular
locus as the centers of the blowing ups. For surfaces, the situation already 
becomes significantly more complicated: the singular loci may be 
1-dimensional and may themselves be singular. In contrast to Hironaka's
approach, the original proof of resolution of surface singularities, of which
the ideas date back to Jung's 1908 article (\cite{Jung},\cite{Walker}), 
does not only involve blow-ups, but also a different type of birational map: 
normalization. 

\subsection{Resolution of surfaces by Jung's approach}

The general philosophy of Jung's approach is to consider a projective 
surface as a ramified covering of the projective plane. The branch locus,
i.e. the discriminant curve, of this covering can then be desingularized
until it only has the simplest possible singularities -- that is until it 
is normal crossing and all components are non-singular. Modifying the
ramified covering in this way and then passing to the normalization of 
the surface, the only remaining singularities of the resulting normal
surface lie over the singularities of the resolved discriminant curve
and are toroidal. These singularities can then be resolved combinatorially.\\

In comparison to the projection, the resolution of the branch locus and
the eventual combinatorial resolution of the toroidal singularities, the
most expensive step here is the normalization. \\

\begin{example}
\begin{verbatim}
> LIB "resjung.lib";          // load the library resjung
> ring R=0,(x,y,z),dp;        // ring of char 0 with 
                              // variables x,y and z
> ideal J=x2+y3;              // ideal of the surface 
> list M=jungresolve(J,0);    // compute the resolution
> M;                          // the output is a list of rings  
[1]:                          // here happens to be only one chart
// characteristic : 0
// number of vars : 4
//      block   1 : ordering dp
//                : names    T(1) x y z
//      block   2 : ordering C
\end{verbatim}
\end{example}    

\subsection{Beyond the geometric case:\\ Lipman's construction for two 
          dimensional schemes}

Like Jung's approach, Lipman's approach makes use of normalization.
But here the surface is not prepared before normalizing it, instead 
normalization steps may occur at many different points in the algorithm.
More precisely, given an excellent, noetherian, reduced 2-dimensional 
scheme, this algorithm proceeds by first normalizing, then blowing up
in the singular locus and then continuing in this way by normalizing 
whenever the resulting scheme at an intermediate step is not normal
and by blowing up the result of this normalization in its singular locus
again.\\
 
As normalization is a more expensive computation than a blow up, this
can potentially slow down the algorithm significantly in comparison to
Jung's algorithm. On the other hand, this approach provides the advantage 
of being applicable for a wider range than the other one 
even opening up the possibilities for considering arithmetic surfaces.\\

\begin{example}

\begin{verbatim}
> LIB "reslipman.lib";                        // load the library reslipman
> ring r=0,(x(1..3)),dp;                      // ring of char 0 with the
                                              // variables x(1),x(2),x(3)
> ideal J=x(2)^2*x(3)^2-x(1)^3;               // ideal of the surface
> list rl=lipmanresolve(J);                   // compute the resolution
> size(rl);                                   // there are 3 charts
3       
\end{verbatim}

If we had done the calculation by hand, we would have used one normalization
and one subsequent blowing up whose center is the only singular point of the
normal surface.\\
As before, we may pass to any of the charts and have a look at the data 
there:

\begin{verbatim}
> def S=rl[1];     // first chart
> setring S ;      // set this ring  
> slocus(J_new);   // compute the singular locus
_[1]=1             // non singular
\end{verbatim}
\end{example} 

\section{Comparisons and Timings} \label{sec5}

The following table shows several examples computed with resbinomial, resolve, resjung and reslipman, over fields of characteristic zero. The computation times $t$ are measured in seconds, and the $\emph{size}$ means the size of the output in terms of number of final affine charts. The subscripts $b$, $s$, $j$ and $l$ indicates resbinomial, resolve, resjung and reslipman respectively.

\begin{table}[h] 
\caption{Test examples}
\begin{center} 
\begin{tabular}{||c|c|c|c|c|c|c|c|c||} \hline
surface or variety & $t_{b}$  &  $t_{s}$ &  $t_{j}$ &  $t_{l}$ & $size_{b}$  &  $size_{s}$ &  $size_{j}$ &  $size_{l}$ \\ \hline
$x_1^2-x_2x_3^2$ & $0$ & $0$ & $1$ & 0 & $12$ & $1$ & $3$ & 1 \\ \hline & & & & & & & & \vspace*{-0.3cm} \\ \hline
$x_1^2-x_2^2x_3^2$ & $0$ & $1$ & $1$ & 0 & $7 $ & $4$ & $2$ & 1 \\ \hline
& & & & & & & & \vspace*{-0.3cm} \\ \hline
$x_1^2+x_2^3x_3^3+x_2^2x_3^5$ & $-$ & $1 $ & $1$ & 1  & $-$ & $14$ & $5$ & $3$ \\ \hline
& & & & & & & & \vspace*{-0.3cm} \\ \hline
$x_1x_2^2x_3^3-x_4^6$ & $9$ & $>7200$ & $-$ & $-$  & $240$ &  & $-$ & $-$ \\ \hline
$x_1^2-x_2^2x_3^2, x_4^2+x_2^3 $ & $1$ & $2$ & $>7200$  & 1   & $43$ & $10$ &  &  1 \\ \hline
& & & &  & & & & \vspace*{-0.3cm} \\ \hline
$x_1^3-x_2x_3^3,x_2^4+x_1^2x_4^2 $ & $313$ & $>7200$ & $>7200$  & 3 & $1863$ &  &  & 4 \\ \hline
\end{tabular}
\end{center} 
\begin{center} 
$t$= time in seconds. \emph{size}=number of final charts.
\end{center}
\end{table} 

Test results with hardware: Intel core i3 processor 330M (2.13 GHz, 1066 MHz FSB),
4 GByte memory. \\

%\vspace*{1cm}

The computed test examples show the stability of resbinomial with respect to the equations generating the ideal. Since the algorithm is basically a combinatorial game on the exponents, similar equations provide similar size of the output-- increasing with growing size of exponents. Note that the number of charts changes considerably also when passing to higher dimensional varieties.\\

In the case of resolve we can see the influence of the geometry of the problem. It is more efficient than resbinomial if at the beginning is possible to consider a permissible center which is bigger than the center selected by resbinomial. \\

%In the rest of cases, for binomial ideals, Villamayor algorithm is less efficient than the specific algorithm for binomial ideals.  \\

With respect to resjung, the number of charts remains almost constant for the examples of hypersurfaces in ambient dimension three. \\

In Lipman's approach, $size_{l}=1$ means that we finish after the first normalization. This approach is the most efficient in the computed examples in ambient dimension $3$, but this is mainly due to the rather simple structure of the chosen examples, which are very accessible to Grauert's normalization algorithm, the bottleneck of Lipman's approach. \\

\end{document}